\documentclass[smallextended]{svjour3}
\usepackage{amsfonts}
\usepackage{latexsym}
\usepackage{amsmath}
\usepackage{amssymb}
\usepackage{amscd}
\usepackage{epsfig}
\usepackage{graphicx, subfigure}
\usepackage{colortbl}
\addtocounter{MaxMatrixCols}{4}
\renewenvironment{proof}{\par {\sc {\bf Proof.}\hskip 5pt}}{\hfill \qed \par}
\usepackage{color}
\usepackage{verbatim}

\usepackage{float}
\floatstyle{plain}
\newfloat{algorithmic}{htbp}{loa}[section]
\floatname{algorithmic}{Algorithm}

\allowdisplaybreaks[1]

\def\nthinsp{\mskip -2   mu}
\def\superstar{^{\raise 0.5pt\hbox{$\nthinsp *$}}}
\def\xstar{x\superstar}
\def\minimize#1{{\displaystyle\minim_{#1}}}
\def\minim{\mathop{\hbox{\rm minimize}}}
\def\subject{\hbox{\rm subject to}}
\def\mtx#1#2{\renewcommand{\arraystretch}{1.2}%
      \left(\! \begin{array}{#1}#2\end{array}\! \right)}
\def\Ahat{\widehat A}
\def\T{^T\!}

\def\A{_{\scriptscriptstyle A}}

\def\xstar{x^{\raise 0.5pt\hbox{$\nthinsp *$}}}

\title {Finding a Hamiltonian cycle by finding the global minimizer of a linearly constrained problem}
	
\author{M. Haythorpe and W. Murray\thanks{The material contained in this manuscript is based
		upon research supported by the
		National Science Foundation Grant DDM-8715153
		and the Office of Naval Research Grant N00014-90-J-1242.}}
\institute{M. Haythorpe
\at Flinders University, Adelaide, Australia\\
\email{michael.haythorpe@flinders.edu.au}
\and
W. Murray
\at Stanford University, Palo Alto, US\\
\email{walter@stanford.edu}
}

\date {}

\begin{document}


%

\thispagestyle{plain}

\maketitle \thispagestyle{empty}

\begin{abstract}

It has been shown that a global minimizer of a smooth determinant of a matrix function corresponds to the largest cycle of a graph.
When it exists, this is a Hamiltonian cycle. Finding global minimizers even of a smooth function is a challenge.
The difficulty is often exacerbated by the existence of many global minimizers. One may think this would help, but in the case of Hamiltonian cycles the ratio of the number of global minimizers to the number of local minimizers is typically astronomically small. There are various equivalent forms of the problem and here we report on two.

Although the focus is on finding Hamiltonian cycles, and this has an interest in and of itself, this is just a proxy for a class of problems that have discrete variables. The solution of relaxations of these problems is typically at a degenerate vertex, and in the neighborhood of the solution the Hessian is indefinite. The form of the Hamiltonian cycle problem we address has the virtue of being an ideal test problem for algorithms designed for discrete nonlinear problems in general. It is easy to generate problems of varying size and varying character, and they have the advantage of being able to determine if a global solution has been found. A feature of many discrete problems is that there are many solutions. For example, in the frequency assignment problem any permutation of a solution is also a solution. A consequence is that a common characteristic of the relaxed problems is that they have large numbers of global minimizers and even larger numbers of both local minimizers, and saddle points whose reduced Hessian has only a single negative eigenvalue. Efficient algorithms that seek to find global minimizers for this type of problem are described. Results using BONMIN, a solver for nonlinear problems with continuous and discrete variables, are also included.

  \smallskip

  Keywords: Global optimization, discrete optimization, binary variables, negative curvature, barrier functions, Hamiltonian cycles.
\end{abstract}

\section{Introduction}
Given an undirected graph $\Gamma$ containing $N$ nodes, determining whether any simple cycles of length $N$ exist in the graph solves the Hamiltonian cycle problem (HCP).
Simple cycles of length $N$ are known as Hamiltonian cycles (HCs). The Hamiltonian cycle problem, named for Sir William Rowan Hamilton's creation of the {\em Icosian} game based on finding Hamiltonian cycles in the planar projection of a dodecahedron, was first formally studied by Thomas Kirkman in the 1800s, although Euler's consideration of the Knight's Tour problem a century earlier can be viewed in the context of HCP. In more recent times, HCP has gained attention both because it was among the first problems proved to be NP-complete \cite{karp}, and also because of its close relationship to the arguably more-famous Traveling salesman problem. The latter can be thought of as the problem of trying to find an optimal Hamiltonian cycle in a weighted graph. Algorithms designed to solve HCP have tended to take the form of either edge-exchange heuristics (e.g. Helsgaun \cite{LKH} or Baniasadi et al. \cite{slh}), branch and cut techniques (e.g. Applegate et al. \cite{concorde} or Chalaturnyk \cite{chalaturnyk}), or genetic algorithms (e.g. Nagata and Kobayashi \cite{nagata}). However, approaches based on directly solving nonconvex optimization formulations have, to date, not been widely considered. In this paper, we describe such an approach that attempts to find a Hamiltonian cycle of a graph by finding a global minimizer of a smooth function.
First, we note that although $\Gamma$ is undirected, it is convenient to think of each edge as a pair of directed arcs. Then, we associate a variable $x_{ij}$ with each (directed) arc $(i,j) \in \Gamma$. Define a matrix $P(x)$, whose $(i,j)$th element is $x_{ij}$ if $(i,j) \in \Gamma$, or 0 otherwise.

It was shown by Ejov et al. \cite{detpaper} that a longest cycle of a graph corresponds to a global minimizer of the problem:
\begin{equation}                                      \label{eqn-det}
 \begin{array}{l@{\hspace{10pt}}r@{\hspace{4pt}}c@{\hspace{4pt}}l}
\minimize{x} & f(x) \equiv -\det(I - P(x) + \frac{1}{N} e e^T) & &                     \\[1.5ex]
	\subject     &   P(x) \in \mathcal{S} , \quad x \ge 0,
 \end{array}
\end{equation}
where $e$ is an appropriately-sized column vector full of units, and $\mathcal{S}$ is the set of stochastic matrices, that is, matrices with nonnegative entries and unit row sums. We shall refer to the linear constraints that arise from this restriction on $P(x)$ as the  $\mathcal{S}$ constraints. It follows that we may also restrict $P(x) \in \mathcal{DS}$, where $\mathcal{DS}$ is the set of doubly stochastic matrices (in which column sums are also one), since if a HC exists and the solution to (\ref{eqn-det}) is defined to be $\xstar$ then $P(\xstar)$ is a permutation matrix. It is these two forms of the problem that we investigate. Note that constraints of these forms are common in other types of discrete problems. The elements of $P(\xstar)$ that are 1 denote the arcs in the HC. We have introduced $x$ as the nonzero elements of $P(x)$. In practice it is convenient to sometimes refer to $x$ as an element of a matrix (ie $x_{ij}$) and sometimes as the vector $x$ composed of those elements. For example, we introduce the concept of twin variables, specifically $x_{ij}$ is the twin of $x_{ji}$. There is no simple way of doing that when $x$ is a vector. Conversely, since we define $f(x)$ in an optimisation problem, $x$. We number the indices by row. So if there are 3 arcs from the first node, there will be $x_1$, $x_2$ and $x_3$ in the first row  of $P(x)$. The first nonzero element of $P(x)$ in the next row is $x_4$ and so on. $P(x)$ is symmetric in the pattern of the nonzero elements but is not a symmetric matrix. The elements in the upper triangular half correspond to arcs in one direction, and their lower triangular half reflection is the arc being taken in the opposite direction.

In what follows, for the sake of simplicity of notation, we will drop the dependency on $x$ and simply refer to $P(x)$ as $P$. We define the term ``leading principal submatrix" as the $(n-1) \times (n-1)$ submatrix obtained by deleting the last row and column of a $n \times n$ matrix. When $P \in \mbox{Int}(\mathcal{DS})$, it is possible to replace the objective function in (\ref{eqn-det}) by the negative determinant of the leading principal submatrix of $I - P$. The result follows from the fact that the restrictions we place on $x$ and $P$ ensure that the LU factorization of $I - P$ exists without the need to permute the rows or columns. Also $I - P$ has rank $N - 1$ and the leading principal submatrix is full rank. The proof was given by Filar et al. \cite{FHM}.

Using the leading principal submatrix has the advantage that the rank-one modification
$\frac{1}{N} e e^T$ is not required, which makes calculating the gradient and the Hessian simpler. A method for efficiently computing the gradient and Hessian of the negative determinant of the leading principal submatrix was provided by Filar et al. \cite{FHM}, and was proved by Haythorpe \cite{michaelthesis} to be more numerically stable than the objective function in (\ref{eqn-det}). Another benefit is that the maximum value is independent of the size of the graph, eliminating the need to scale any parameters by the size of the graph.


In the $\mathcal{DS}$ case the problem of interest is of the form
$$\minimize{x} \;\;\;\; f(x)$$
subject to
\begin{eqnarray}P & \in & \mathcal{S},\label{eq-DS1}\\
P^T & \in & \mathcal{S},\label{eq-DS2}\\
x_{i} & \geq & 0, \quad\forall (i),\label{eq-DS3}\end{eqnarray}

where $f(x)$ here is the negative determinant of the leading principle submatrix of $I - P$.
Constraints (\ref{eq-DS1})--(\ref{eq-DS3}) are called the {\em doubly-stochastic} constraints. For brevity, we refer to constraints (\ref{eq-DS1})--(\ref{eq-DS3}) as the $\mathcal{DS}$ constraints. It is assumed that any graph considered is simple and undirected.
Although this is a classical linearly constrained problem, it is different in character from those whose variables are not related to a binary-variable problem. Discrete problems that are transformed to that of finding a global minimizer in a continuous space have common characteristics. First note that if there is only one global minimizer and no local minimizers the problems would be trivial to solve. Since they are not trivial to solve, these problems have many minimizers, including potentially many global minimizers. One consequence of having multiple minimizers is the problem is likely to have many saddle points. A particular concern is that it can be shown that for any two isolated minimizers, there exists a path between them that has a stationary point with a Hessian with only one negative eigenvalue. In other words, this point is very similar to a minimizer. A proof of this result is provided in the supplementary material for this paper. Moreover, as the number of isolated minimizers increases there is the potential for the number of such points to increase quadratically in the number of global minimizers. In the Hamiltonian cycle problem it is known that in general,  the number of solutions grows exponentially with $n$, with each corresponding to a unique global minimizer \cite{eppstein}. While this may make such a transformation of a discrete problem seem unattractive, there is an advantage over algorithms finding solutions whose variables are not discrete, namely we can round to the discrete choices at any iteration and check whether there is a solution. This avoids the need to iterate to a limit, along with the conditioning problems that can sometimes cause.

It was shown by Filar et al. \cite{FHM} how to compute $f(x)$ and its first and second derivatives very efficiently. This is critical since we show that directions of negative curvature are essential to solving this problem, and they play a much more critical role than is typically the case. A key issue is symmetry. Obviously, for \emph{every} HC there is another HC obtained by reversing the direction. This symmetry reveals itself in the problem variables. If there is a HC with $\xstar_1 = 1$ then there exists a reverse cycle that is also a HC, in which $\xstar_1 = 0$ and its twin variable is 1. A critical characteristic of the approach we advocate is to avoid introducing bias in the initial estimate of the solution. Consequently, we set the initial value of these two variables to be identical (they may both be 0 in another HC). Quite frequently (and this almost always happens with some pairs), when using only descent, some of these twin variables remain equal. In such circumstances it is only the use of a direction of negative curvature that breaks the tie. While such behavior is possible for nonlinear optimisation problems, it is quite rare. Consequently, in this class of problem directions of negative curvature play a more important role, and often more important than that of using a direction of descent. Again, unlike standard nonlinear optimisation problems where we usually observe no directions of negative curvature in the neighborhood of the solution, here they are always present, which is one reason why the solution is at a vertex. What is happening is that from our current iterate there are two equally attractive minimizers, so it steers a course going to neither unless directions of negative curvature are used. The symmetry reveals itself also in the problem function and derivatives. If the iterates to solve the problem are denoted by $\{x_k\}$, then at $x_k$ it is usually the case that $g_k$, the gradient of $f(x)$ at $x_k$, is orthogonal to the eigenvectors corresponding to negative eigenvalues of the Hessian of $f(x)$ at $x_k$. Consequently, when solving for the Newton step using the conjugate gradient algorithm, it will not detect when the Hessian is indefinite. Again, for standard nonlinear optimisation problems this is extremely rare.

The primary contributions of this manuscript are as follows. First, we extend the result that the LU factorization for $I - P$ exists without requiring permutations for $P \in Int(\mathcal{DS})$, to the more general case where $P \in \mathcal{S}$. Second, we present an interior point algorithm designed to find Hamiltonian cycles, and analyze its performance. This analysis reveals a number of unusual features of this kind of problem. We propose that these features should be taken into account when designing optimization algorithms for discrete problems that relax the discrete requirement.

The remainder of this manuscript is organized as follows. In section 2 we consider further the two different choices of constraints that may be used. In section 3 we discuss different equivalent forms of the problem, in particular considering the determinant of leading principle submatrix of $I - P$, and investigate the nature of the LU factorization of this matrix. In section 4 we describe the approach we use to solve this problem via an interior point method. In section 5 we give detailed descriptions of the various algorithms that we use to solve the problem. In section 6 we provide some numerical results, including a comparison with BONMIN \cite{bonmin1,bonmin2}, a powerful nonconvex solver. Finally, section 7 contains some concluding remarks.

\section{Choice of constraints}

Typically, in an optimization problem it is better to have more constraints if such constraints can be added, even if these are inequalities and are known not to be active at the solution. However, it is not always the case, as there can be a trade-off between the benefits of reducing the search space, and the inefficiency of making the iterates more complicated to compute. In this problem, we have a choice of either $P \in \mathcal{DS}$ or $P \in \mathcal{S}$ (obtained by either eliminating equations (\ref{eq-DS1}) or (\ref{eq-DS2})). Note here that when solving with $P \in \mathcal{DS}$ we are adding more equality constraints without adding extra variables, and hence we are reducing the degrees of freedom in the problem. It is therefore worth taking time to consider which set of constraints is better to work with.

It was shown by Ejov et al. \cite{detpaper} that when $P \in \mathcal{DS}$, the LU factors of $I - P$ exist regardless of the pivoting order. This has many beneficial consequences, not the least of which is the objective of the problem may be recast to be $-\det(M)$, where $M$ is the leading principle submatrix of $I - P$. Although $I-P$ is singular, its leading principle submatrix is nonsingular as a consequence of the existence of the LU factorization. Note that $I - P$ is typically very sparse (if it is not then the underlying graph is dense, and finding a HC is usually trivial). However, we shall show in Section \ref{sec-prelim} that a similar property is also true even when $P \in \mathcal{S}$, so this is not a reason for preferring  $P \in \mathcal{DS}$. It was shown by Ejov et al. \cite{detpaper} that when $P \in \mathcal{S}$,  if a variable in $P$ is not 0 or 1, then altering it to one of them reduces the objective. One consequence is that all local minimizers are binary. Also, a common practice when solving relaxed problems is rounding. For unconstrained problems this is trivial, but for constrained problems, rounding may cause infeasibility. However, with $P \in \mathcal{S}$, it is clear that rounding is trivial even though the problem is constrained.

When $P \in \mathcal{DS}$, changing a variable to either 0 or 1 does not necessarily improve the solution. However, it will be seen that one of the steps we propose in our algorithm is deletion or deflation, which occur when one or more of the variables is set to 0 or 1 respectively. In order to do this for $P \in \mathcal{DS}$ without losing feasibility, either an linear program (LP) or quadratic program (QP) needs to be solved to obtain an appropriate feasible value. Since a strictly interior point is required, and sometimes one does not exist, there are some complications involved with this task.

Our algorithm requires finding both a sufficient feasible descent direction and a feasible direction of sufficient negative curvature. To ensure feasibility, we define a null-space matrix $Z$, such that $A Z = 0$, where $A$ is the matrix of constraint coefficients. The matrix $(A Z)$ is full rank. Typically, such $Z$ matrices are almost always dense.  Consequently,  the smaller the dimension of $Z$ the better. However, it was shown by Filar et al. \cite{FHM} that there exists a $Z$ for  the $\mathcal{DS}$ case that is sparse and structured. A sparse $Z$ can also be constructed for the $\mathcal{S}$ case that, despite being larger than for the $\mathcal{DS}$ case, is even simpler and sparser. Note that linear equality constraints of this character are common to other discrete problems.

Given that it is unclear, from a preliminary analysis, which set of constraints is better to use, we will begin by considering both. As will be seen from the upcoming experimental results, the $P \in \mathcal{DS}$ constraints prove to be superior despite some of the additional complexities involved with using them.

\section{Preliminaries}\label{sec-prelim}

As already noted, it is beneficial if the determinant in the $\mathcal{S}$ case may be computed from a LU factorization only of $I - P$, since for most difficult problems this is sparse. For large problems the known fixed pattern of nonzeros in the factorization enables the factors to be found more efficiently than would be the case if numerical pivoting was required.

\subsection{The LU factorization $I - P$}
We show that an LU factorization of $I - P$ and of $I - P^T$ exists when $P$ is a stochastic matrix, such that the last diagonal element of the resulting upper factorization is zero. As already noted this was shown to be true for a doubly stochastic matrix. Since $I - P$ is singular, this result implies that we can dispense with the term $(\frac{1}{N} e e^T)$, and require only the determinant of the sparse leading principle submatrix of $I - P$.

A stochastic matrix may have either rows or columns that sum to unity. In forming the LU factorization, it is common to assume row interchanges rather than column interchanges. This is just convention and there is no advantage to doing it one way or the other. However, for sparse matrices the manner in which the few elements are stored does matter when performing the LU factorization. Since, when forming such matrices, it is assumed that row interchanges will be done, this has an impact on how best to store the sparse matrix in compact form. Note also that the original problem could be posed with an objective with replaced by $I - P^T$. In the proof we assume row interchanges may be made, and this causes us to prefer to assume that $P$ has unit columns. The converse result for unit rows follows immediately from this result.

\begin{definition} A  $n\times n$ matrix $A$ is said to have
property $S_c$ if
\begin{enumerate}
\item[] $A_{i,i} \ge 0$ for $\forall$ $i$,

\item[] $A_{i,j} \le 0$ for $\forall$ $i \ne j$,

\item[] $ A^T e = 0$.

\end{enumerate}

\end{definition}

\begin{theorem} If $A$ has property $S_c$, an LU factorization of $A$ exists such that $U_{n,n} = 0$.
\end{theorem}
\begin{proof} We provide an inductive argument that shows if $A$ has property $S_c$, then each subsequent step of Gaussian Elimination produces a new submatrix which also has property $S_c$ to be considered in the next iteration.

First, consider the case when $A_{1,1} = 0$. By definition, in this case the whole first column is zero. Then, there is an option for Gaussian elimination (GE) on how to proceed. Specifically, we have the option of choosing what portion of the first row we add to the remaining rows. The rule we follow here is that the sum of the proportions (all non-negative) that we add should be 1, and the simplest strategy for achieving this is to just add the first row to the second. This rule should then followed throughout should subsequent pivots be zero. It can be seen that when this occurs, then the column of $L$ generated at the $i$th step is $e_i^T - e_{i+1}^T$ and it follows that $e^T ( e_i - e_{i+1}^T) = 0$.

Note that $A_{1,1}$ is an element of largest magnitude in the first column of $A$. In general, after one step of Gaussian elimination (GE) we get
$$
     A = L_1 \mtx{cc}{1 & 0 \\ 0&\Ahat}\mtx{cc}{U_{1,1} & u^T \\ 0 &I},
$$
where
$$
      L_1 = \mtx{cc}{1 & 0 \\ l &I}.
$$
Note that $l = -e_2$ in the case when $A_{1,1} = 0$. Now, we have $A^T e = 0$, which implies that
$$
    e\T \mtx{cc}{1 & 0 \\ l &I}\mtx{cc}{1 & 0 \\ 0&\Ahat}\mtx{cc}{U_{1,1} & u^T \\ 0 &I} = 0.
$$
It follows that $\Ahat^T e = 0$ and $e\T l = -1$. By definition, we have
$$
    \Ahat_{i,j} = A_{i+1, j+1} - l_iu_j \quad \forall \quad i\ne j.
$$
Since $l_i \le 0$ and $u_j \le 0$ it follows that $\Ahat_{i,j} \le
0$ $\forall$ $i \ne j$. From this result and $\Ahat^T e = 0$ it
follows that $\Ahat_{i,i} \ge 0$ and that $\Ahat$ has property
$S_c$. We can now proceed with the next step of GE.

In the final iteration we have that $\Ahat $ is the single element $U_{n,n}$, and still has property $S_c$. Hence, it follows that $U_{n,n} = 0$
\end{proof}

\begin{corollary}
$L^T e = e_n$.
\end{corollary}

\begin{corollary}
When performing GE, the elements being eliminated are not larger in magnitude than the pivot. This implies that
$0 \ge L_{i,j} \ge -1$ for $i \ne j$.
 \end{corollary}

 This property implies that $L$ is about as well conditioned as it can be. Moreover, if standard software to perform GE is used, even if it performs row interchanges when needed, they will never be required and the LU factorization  of $A$ will be obtained and not that of $\Pi A$, where $\Pi$ is a permutation matrix. At each iteration of our algorithm the fill-in of the LU factors is identical.


\begin{definition} A $n\times n$  $A$ is said to have
property $S_r$ if
\begin{enumerate}
\item[] $A_{i,i} \ge 0$ for $\forall$ $i$

\item[] $A_{i,j} \le 0$ for $\forall$ $i \ne j$

\item[] $ A e = 0$.

\end{enumerate}

\end{definition}

\begin{corollary}
If $A$ has property $S_r$ then an LU factorization of $A$ exists.
\end{corollary}
Note the leading principle submatrix of $A$ is the transpose of that of $A^T$.

\begin{lemma}
  If $A$ has no off diagonal entries that are zero, then A is rank $n - 1$ and the only diagonal of $U$ that is zero is $U_{n,n}$.
\end{lemma}
\begin{proof} Assume $n \ge 3$. Note that $A_{1,1} > 0$ and is an element of largest magnitude in the first column of $A$. After one step of Gaussian elimination (GE) we get
$$
     A = L_1 \mtx{cc}{1 & 0 \\ 0&\Ahat}\mtx{cc}{U_{1,1} & u^T \\ 0 &I},
$$
where
$$
      L_1 = \mtx{cc}{1 & 0 \\ l &I}.
$$
 By definition, we have
$$
    \Ahat_{i,j} = A_{i+1, j+1} - l_iu_j \quad \forall \quad i\ne j.
$$
Since $l_i < 0$ and $u_j < 0$ it follows that $\Ahat_{i,j} <
0$ $\forall$ $i \ne j$. We have that $\Ahat^T e = 0$, which implies that
$\Ahat_{1,1} > 0$. Inductively, we can repeatedly apply the same reasoning until the dimension of $\Ahat$ is 1.\end{proof}

\begin{corollary}
If $A$ is rank $n-1$ then the leading principle submatrix is nonsingular.
\end{corollary}

In reality, the matrix of interest in this manuscript, $A \equiv I -P$, has stronger properties than  property $S_c$. For example its diagonal elements are all units, the non-zeros form a symmetric pattern and for the iterates generated in the algorithm its rank is $n - 1$ since the determinant is never zero. Using these last two properties we can show that det($\Ahat$), where $\Ahat$ is the leading principle submatrix of $I - P$ is a surrogate for our original objective function.

\begin{lemma} If $A \ne 0$ is symmetric in its pattern of zero elements and A is rank $n - 1$ then the leading principal submatrix of $A$ is nonsingular.
\end{lemma}
\begin{proof}
First, assume that $A_{1,1} = 0$. Then all the elements in the first column are 0 and by the symmetric assumption so are the elements in the first row. However, since the  $\Ahat e = 0$ this implies the rank of $A$ is less than $n-1$, so this assumption is false. We can repeat this process inductively to show that $\Ahat_{1,1} \ne 0$ until the dimension of $\Ahat$ is 1. Therefore, the only diagonal of $U$ that is zero is $U_{n,n}$, and hence the result follows.\end{proof}

Another useful feature of the matrices we are interested in is that the zero element pattern of $L$ is always the zero element pattern of the lower part of $A$ plus some additional elements. It cannot be the case that there is a zero element of $L$ in the location of a nonzero of $A$ resulting from the actual numerical values of $A$, which can happen for general matrices. The nonzeros in $U$ are in the same location as the nonzeros in $L^T$. This can be useful when computing the LU factors and saves storing two set of indices.

\subsection{Equivalent objective functions}
Consider the following objective function
 \begin{equation} \label{problemee}
  \det(A + \alpha ee^T),
\end{equation}
where $A = I - P$ and $\alpha > 0$. It can be shown that if the
eigenvalues of $A$ are $\lambda_1, \dots, \lambda_{n-1}, 0$, then the eigenvalues of $A + \alpha ee^T$ are $\lambda_1, \dots, \lambda_{n-1}, n\alpha$.

An equivalent objective function to (\ref{problemee}) is:
 \begin{equation} \label{problemeen}
 \det(A + \alpha ee_n^T).
\end{equation}
 The eigenvalues of  $A + \alpha ee_n^T$  are
$\lambda_1, \dots, \lambda_{n-1}, \alpha$. Since the
eigenvalue that differs is independent of $x$, the problems are
equivalent. Note that it follows that $A + \alpha ee_n^T$ is
nonsingular when $A$ is rank $n-1$.

This result is important since it implies that, if the LU factorization of $A$  exists without pivoting {\em and} the last diagonal element of $U$ is zero, then
\begin{equation}
       \det(\Ahat),
\end{equation}
where $\Ahat$ is the leading principle submatrix of $A$, is also an equivalent objective function to (\ref{problemee}) and (\ref{problemeen}) by the following argument. First note that $\det(\Ahat) = \det(\bar A)$, where
$$
        \bar A =  \mtx{cc}{\Ahat & 0 \\ 0 &1}.
$$
Let $LU$ denote the LU decomposition of $A$, and $\bar L \bar U$ denote the LU decomposition of $\bar A$. Since, in the upper triangle, only the last column of $A$ differs from the columns of $\bar A$ it follows the only diagonal element of U that is different in the two cases is $U_{n,n}$ and $\hat U_{n,n}$. The part of the elements in last column of $\Ahat$ that originates from the last column of $A$ contribute zero to the value of $\hat U_{n,n}$. It follows that $\det(\Ahat) = \gamma \det(A + \alpha ee_n^T)$, where $\gamma$ is a positive constant.

\section{Finding the global minimizer of the linearly constrained problem}

The basic approach used is similar to that due to Murray and Ng \cite{murrayng}, who first relax the problem, and then solve a sequence of problems in which a strictly convex function is added to the objective together with a nonconvex function that attempts to force the variables to be binary. Initially, the strictly convex function dominates the objective, and in the limit the nonconvex term dominates. Our approach is a simplification since the objective function we use is minimized at binary points, so the nonconvex term is not needed. Applying this general approach to a specific problem with significant structure, the algorithm can be modified to improve not only efficiency, but also to improve the likelihood of obtaining a global minimizer and hence a HC.

Our key focus is on how the individual problems in the sequence are solved. We wish to capitilize on the knowledge that the solution is at a degenerate vertex, and the Hessian of the objective is indefinite even in the neighborhood of the solution. It is our thesis that this shifts the relative importance of descent directions and directions of negative curvature. It will be seen that a much heavier use of negative curvature is made, with less emphasis on the use of descent directions; this is the reverse of what optimization algorithms usually do. A peculiarity, which we think may be true of most problems with multiple global minimizers, is that the gradient at the iterates is often spanned by the eigenvectors corresponding to the positive eigenvalues of the Hessian, even though the Hessian is indefinite. This corresponds to the so-called ``hard case" in trust-region methods. Typically, in such methods little or no attention is paid to this situation since it is considered very unlikely to arise, and essentially impossible to keep arising; however, that does not appear to be the case here.

In both the stochastic and doubly stochastic case, we are interested in solving a problem of the form:
\begin{eqnarray}& \min \; f(x) &\nonumber\\
& \mbox{s.t.} &\label{eq-problem0}\\
& Ax = e, & i = 1, 2, \hdots, m,\nonumber\\
& 0 \le x \le e.\nonumber\end{eqnarray}
Strictly speaking, the upper bounds on $x$ are not required since the equality constraints and the lower bounds ensure that the upper bound on $x$ holds. However, for now we shall leave them in.

We are interested in the global minimizer, and a typical descent algorithm will only converge to the local minimizer associated with the initial point. Murray and Ng \cite{murrayng} propose adding a strictly convex function $\mu \phi(x)$ to the objective, where $\mu$ is a positive scalar. A sequence of problems is then solved for a sequence of strictly monotonically decreasing values of $\mu$. For $\phi(x)$ with certain continuity properties, the trajectory of minimizers $\xstar(\mu)$ is unique, continuous, and smooth. When the initial value of $\mu$, say $\mu_0$, is sufficiently large, the new objective is also strictly convex and has a unique and therefore global minimizer. Hence, in such a case, given any initial point we will converge to a unique optimal point $\xstar(\mu_0)$. Then, the trajectory $\xstar(\mu)$ is a smooth continuous arc for $\mu \in [0, \mu_0]$. When the initial choice of $\mu_0$ is not sufficiently large then the solution will depend on the initial point. Every local minimizer of the original problem has an associated arc. The {\em range} of $\mu$ for which these arcs exist differs depending on the local minimizer. By choosing an initial $\mu_0$ large enough, we eliminate the arcs corresponding to local minimizers that are not linked to the global minimizer of the barrier function. In other words, the arc corresponding to the global minimizer is the last arc left standing as $\mu$ is increased. Consequently, a minimizer found by this algorithm is the one whose trajectory is linked to the initial unique global minimizer.

A feature of the problem is that it has what we term ``twin variables".  From the definition of the variables as elements of the matrix $P$,  if $P_{ij}$ is not always zero then neither is its twin $P_{ji}$. In terms of the graph, this corresponds to the same edge except in the opposite direction. Since the reverse of a HC is itself a HC, twin variables have an equal probability of being in a HC irrespective of the structure of the underlying graph. It is essential in the approach that we consider that the minimizer of $\phi(x)$ is a neutral point with regard to the minimizers of the original problem. If it is not the case, then when $\mu$ is large the solution is closer to some local minimizer than another, and we do not know whether this is the global minimizer. For example, consider the problem

\begin{eqnarray}& \min \; f(x) \equiv -(x-0.5e)^T(x-0.5e) + \epsilon c\T x &\nonumber\\
& \mbox{s.t.} &\label{eq-problem}\\
& 0 \le x \le e.\nonumber\end{eqnarray}
When  $\epsilon$ is very small typically this function has local minimizers at all vertices of the feasible region.
The transformed problem is
$$
\min f(x) + \mu \phi(x)
$$
for a sequence $\mu = \mu_0, \mu_1, \dots$ , where $\mu_i \to 0$ as $i \to \infty$.

 If choose $\phi(x)$ with a minimizer that is closer to one of the vertices than the others, then typically we will converge to that minimizer as $\mu$ becomes smaller. The neutral point is $0.5 e$. A suitable choice for $\phi(x)$ is

$$
 \phi(x) = -\sum_{i=1}^n  \ln{x_i} + \ln{(1-x_i)}.
$$

We will use $\phi(x)$ defined in this way as our choice of barrier function throughout the remainder of the manuscript. Hence, we have transformed the original problem into minimizing a sequence of barrier functions, which have been used for many years to eliminate inequality constraints from problems. The reason here for using such functions is not eliminating inequality constraints, that is simply a side benefit. Indeed, retaining the inequalities and solving the original problem using, for example, an active set method, would still have been efficient, especially since we do not expect the size of the problems to be extremely large (100,000 variables or more). Our motivation for using barrier functions is to find a global minimizer, and consequently this impacts how the initial $\mu$ is chosen and how it is subsequently adjusted. It also impacts how best to solve the subproblems, since we know the global solution is at a degenerate vertex. Since we seek an initial neutral point, we set $\mu_0 = \infty$ (equivalent to dropping $f(x)$ from the objective).  A test of whether the choice of $\phi(x)$ leads to a neutral initial point when solving our problems is whether the twin variables have the same initial value, and this is observed in the numerical testing of the barrier function when we use it in our numerical testing. Note that although $\xstar_i$ is either 0 or 1, it is not the case that 0.5 is a neutral value in the presence of linear equality constraints.

Since the barrier function removes the need for the inequality constraints, the algorithms requires the solution of a sequence of linearly \textit{equality} constrained optimization problems.

\subsection{Solution of the linearly equality constrained subproblem}

 The choice of method is dictated by the need to converge to points that at least satisfy second-order necessary conditions. This requires the algorithm to determine whether the reduced Hessian is positive semidefinite. To obtain the reduced Hessian matrix, we need the null space matrix $Z$, which is such that $AZ = 0$ and $A^TZ$ is full rank. The reduced Hessian is then given by $Z\T H Z$, where $H$ is the Hessian of $f(x)$.

 We use a line search method based on determining a descent direction, and when available, a direction of negative curvature. A sequence $\{x_k\}$ of improving estimates is generated from an initial feasible estimate $x_0$, from

$$
        x_{k+1} = x_k + \alpha_k (p_k + d_k),
$$
where $\alpha_k$ is a step length that ensures a sufficient decrease, $p_k$ is
a sufficient descent direction,  and $d_k$ is a direction of sufficient negative curvature. It was shown by Forsgren and Murray \cite{forsgrenmurray} that this sequence converges to a point that satisfies the second-order necessary conditions.

Typically, such methods combine a direction of descent with a direction of negative curvature when the latter exists. Our observation is when a direction of negative curvature does exist and is used purely as the search direction, then at almost every subsequent iteration a direction of negative curvature exists and is usually getting stronger. Consequently, when we get a direction of negative curvature, we do not bother computing the direction of descent.

Given the importance of the direction of negative curvature, we depart from normal practice and apply the modified Cholesky algorithm \cite{practical} to the following matrix
$$
    Z\T H Z + \delta I,
$$
where $-\delta$ is an estimate of the smallest eigenvalue of $Z\T H Z$ when it is thought $Z\T H Z$ is indefinite, otherwise $\delta$ is positive and very small in magnitude. The rationale is that when
$Z\T H Z$ is  indefinite, this leads to a very good direction of sufficient negative curvature. When $Z\T H Z$ is positive definite, the small shift ensures that the matrix has a condition number that is sufficiently small to ensure sufficient accuracy in the direction of descent. If a modification is made in the modified Cholesky factorization, then $Z\T H Z + \delta I$ is indefinite and the following system is solved:
$$
     Rd_z = e_j,
$$
where $R$ is the upper triangular factor, and the index $j$ is obtained during the modified Cholesky factorization. It can be shown that $d$,
where $d = Zd_z$ is a direction of sufficient negative curvature. Moreover, we have
$$
  d\T H d \le -\delta d^T d.
$$
We can improve this direction of negative curvature by applying an algorithm minimizing $d\T H d /d^T d$. Specifically, we reduce the value by doing a sweep of univariate minimization of this function. This cost is roughly the same as a matrix-vector multiplication and so can be repeated if need be. We use the improved value as the estimate of $\delta$ in the following iteration. Note that the sign of $d$ is always chosen so that $d\T g \le 0$, where $g$ is the gradient of $f(x)$.

If no modification is made in the modified Cholesky factorization, then $Z\T H Z + \delta I$ is positive definite. If $\delta$ is not small, we cannot say for certain that no negative curvature exists. However, we will know that the smallest eigenvalue is bigger than $-\delta$. We repeat the modified Cholesky factorization with $\delta \leftarrow \delta/2$. If after a small number of reductions we still get no modification, then we assume that $Z\T H Z$ is positive definite, and so we set $\delta$ to the default small value, and compute a direction of sufficient descent by solving
$$
   R^T R p_z = -Z^T g.
$$

We use a very crude line search along either $d$ or $p$. We compute the maximum step to the boundary and take a step $\beta$ times the value. If that is not a lower point, we multiply the step by 0.5 until we succeed. Typically $\beta = 0.9$ and is almost always successful.

\section{The outer algorithm}

A key difference with the use of a barrier function here, compared to solving problems unrelated to relaxed discrete problems, is that $f(x)$ behaves in an unusual way. After a strictly feasible point is found, this is used to minimize the barrier function alone (equivalent to setting $\mu = \infty$). This is an easy function to minimize, and to do so accurately. This is necessary to avoid bias. This differs from what is done when minimizing a regular function, where we would prefer to provide an initial point reasonably close to a solution. Indeed, what we are proposing is to find an initial point as far away as possible from the solutions. In some cases, such as for cubic graphs, the minimizer of the barrier function is known ($x_i = 1/3$). Typically, we want to reduce $\mu$ at a slow rate. However, another feature of the HC problem is that the point that minimizes the barrier function is either a saddle point of the determinant function, or very close to one. Again this rarely if ever happens when using a barrier function for normal problems. The consequence is that moving the iterates from their current location requires reducing $\mu$ sufficiently to make the current reduced Hessian indefinite. Quite how much is not difficult to estimate. The Hessian of the barrier function is a well conditioned diagonal matrix at the minimizer. The condition number is usually less than 2 at this point and for cubic graphs is 1.  In both the stochastic and doubly stochastic case the matrix $Z$ has a low condition number. Consequently, given an estimate of the smallest eigenvalue of either $H_D$, the Hessian of determinant function, or of $Z^T H_D Z$ it is easy to find a good estimate of the change needed in $\mu$. If it is not sufficient then we can simply divide by 10 until it is. In our testing this was never needed. Once we got negative curvature we never reduced $\mu$ again since either we succeeded in finding a HC without needing to, or we failed.

An alternative for solving the standard problem is to use the primal dual approach. The standard approach means that the Newton direction is poor when $x$ is close to a bound. There are two reasons not to use the primal dual method here. Firstly, we do not need to have $\mu$ very small since we know the solution is converging to a binary point, and so we can round and test the solution. Ill-conditioning arises due to a variable becoming close to a bound. Should a variable become close to a bound, that variable can be removed from the problem. How to do this is described in the section on deletion and deflation. Secondly, we need to use directions of negative curvature, but the Hessian in the primal dual formulation is not assured to give a direction of negative curvature except in the neighborhood of a stationary point \cite{dulce}.

\section*{Definition of the null space $Z$}

A common way of defining $Z$, such that $AZ = 0$ and $A^TZ$ is full rank, is to first partition the columns of $A = \left[\begin{array}{cc}B & S\end{array}\right]$, where $B$ is nonsingular. Then we can define

$$Z := \left[\begin{array}{c}-B^{-1}S\\I\end{array}\right].$$


If only stochastic constraints are required, the matrix $A$ can be quite simply defined. If the graph has $N$ vertices, where vertex $i$ has degree $d_i$, then we can define

$$A := \left[\begin{array}{cccc}A_1 & A_2 & \cdots & A_N\end{array}\right],$$

where $A_i = e_i e^T$ is an $N \times d_i$ matrix. It is easy to see that the condition number of $A$ is bounded above by $N$, by checking that $AA^T = \mbox{diag}\left[\begin{array}{cccc}d_1 & d_2 & \cdots & d_N\end{array}\right]$, and therefore the condition number of $A$ is the ratio of the largest degree to the smallest degree.

In order to define the null space matrix for $A$, it is trivial to reorder the columns of $A$ such that the first column of each $A_i$ submatrix appears first. The reordered matrix is $\hat{A} := \left[\begin{array}{cc}I & S\end{array}\right]$, where $S$ is defined as

$$S := \left[\begin{array}{cccc}S_1 & S_2 & \cdots & S_N\end{array}\right],$$

and $S_i = e_i e^T$ is an $N \times (d_i-1)$ matrix. Then the null space matrix for $\hat{A}$ is

$$\hat{Z} = \left[\begin{array}{c}-S\\I\end{array}\right],$$

and appropriately reordering the rows of $\hat{Z}$ provides the null space matrix for $A$. One advantage of defining the null space matrix in this way is that the sparsity inherent in difficult graphs is retained in $Z$, and non-zero entries are all either $+1$ or $-1$. The condition number of $Z$ is equal to $\max_i d_i$. The only operations involving $Z$ that are required are matrix-vector products. Then for a given vector $v$, the product $Zv$ can be computed very efficiently.


If doubly-stochastic constraints are desired, the matrix $A$ defines constraints on both the rows of $P$ and the columns of $P$. We first define a matrix $A_r$, corresponding to the row constraints, to be identical to the $A$ matrix for the stochastic case. Next we define a matrix $A_c$, corresponding to the column constraints. Suppose each variable $x_k$ corresponds to an arc $a_k = (i,j)$. Then $A_c$ is defined as

$$[A_c]_{jk} := \left\{\begin{array}{lcl}1 & & \mbox{if } \exists\; i\;\; \mbox{ s.t. } a_k = (i,j)\\
0 & & \mbox{otherwise}\end{array}\right.$$

Then, we can define $A$ to be

$$A := \left[\begin{array}{c}A_r\\A_c\end{array}\right].$$

The matrix $A$ defined in this way is certain to be rank deficient. In order to construct a null space matrix, we want to delete enough rows to obtain a full rank matrix, and then reorder the matrix to obtain $\hat{A} = \left[\begin{array}{cc}B & S\end{array}\right]$, where $B$ is a \textit{triangular} matrix. This can be achieved by using the following algorithm.

{\scriptsize\begin{center}\begin{tabular}{|l|}\hline
\\
\hspace*{0.4cm}{\bf Input}: $A, N$. \;\;\;\;\; {\bf Output}: $B, S, \mathcal{I}$.\\
\\
\hspace*{0.4cm}{\bf begin}\\
\hspace*{1.1cm}count $\leftarrow$ 0\\
\hspace*{1.1cm}rows $\leftarrow$ rank$(A)$\\
\hspace*{1.1cm}$\hat{A} \leftarrow A$ with rows removed to make $\hat{A}$ full rank\\
\hspace*{1.1cm}cols $\leftarrow$ columns$(\hat{A})$\\
\hspace*{1.1cm}$r \leftarrow $ rows\\
\hspace*{1.1cm}$\mathcal{I} \leftarrow \{1, \hdots, $cols$\}$\\
\hspace*{1.1cm}{\bf while} $r > 0$\\
\hspace*{1.8cm}C $\leftarrow$ Identify a set of columns $\{c_1, \hdots, c_k\}$ such that $\sum\limits_{i=1}^{r}a_{ic_j} = 1,\;\;\;\forall j = 1, \hdots, k$\;\;\;\;\;\\
\hspace*{2.5cm}and $a_{ic_j}a_{ic_k} = 0,\;\;\;\forall i = 1, \hdots, r,\;\;\; j \neq k$\\
\hspace*{1.8cm}{\bf for} $i$ {\bf from} 1 {\bf to} $k$\\
\hspace*{2.5cm}count $\leftarrow$ count + 1\\
\hspace*{2.5cm}$\mathcal{I} \leftarrow \left[\begin{array}{cccccc}\mathcal{I}_1 & \hdots & \mathcal{I}_{\mbox{count}-1} & \mathcal{I}_i & \mathcal{I}_{\mbox{count}} & \hdots\end{array}\right]$ (Moving $\mathcal{I}_i$ into position count)\\
\hspace*{2.5cm}$\hat{A} \leftarrow \hat{A}(\mathcal{I})$ (Moving column $c_i$ to column count) \\
\hspace*{2.5cm}$\hat{A} \leftarrow$ reorder the rows to get a 1 in positive $(r - i + 1,$count)\\
\hspace*{1.8cm}{\bf end}\\
\hspace*{1.8cm}$r \leftarrow $ rows - count\\
\hspace*{1.1cm}{\bf end}\\
\hspace*{1.1cm}$\mathcal{I} \leftarrow $ Reverse the order of the first rows entries in $\mathcal{I}$\\
\hspace*{1.1cm}$\hat{A} \leftarrow \hat{A}(\mathcal{I})$ (Reverse the order of the first rows columns in $A$)\\
\hspace*{1.1cm}$B \leftarrow \hat{A}(1:\mbox{rows},1:\mbox{rows})$\\
\hspace*{1.1cm}$S \leftarrow \hat{A}(1:\mbox{rows},\mbox{rows+1}:\mbox{cols})$\\
\hspace*{0.4cm}{\bf end}\\
\\
\hline\end{tabular}\\
\vspace*{0.25cm}{\bf Reordering $A$ algorithm}\end{center}}

Then the null space for $\hat{A}$ is defined to be

$$\hat{Z} := \left[\begin{array}{c}-B^{-1}S\\I\end{array}\right].$$
Unlike typical problems, the $Z$ constructed in this way is sparse (similar to that of $A$) and does not require the LU factorization of $B$ since $B$ is lower triangular. Moreover, $B$ has elements that are  either 0 or 1. Consequently, operations with $Z$ do not require any multiplication.

\section*{Deletion and deflation}

If, at any stage, one or more of the $x_{ij}$ variables approach their extremal values (0 or 1), we fix these values and remove the variables from the problem. This process takes two forms: {\em deletion} and {\em deflation}, that is, setting $x_{ij}$ to 0 or 1, respectively. Note that we use the term deflation because in practice the process of fixing $x_{ij} := 1$ results in two nodes being combined to become a single node, reducing the total number of nodes in the graph by 1.

Deletion is a simple process of fixing a variable to 0 by simply removing its associated arc from the graph.  When a variable $x_{ij}$ is close to 1, we perform a deflation step by combining nodes $i$ and $j$ together, effectively removing node $i$ from the graph. Then, we redirect any arcs $(k,i)$ that previously went into node $i$ to become $(k,j)$, unless this creates a self-loop arc. After deletion or deflation we construct the new constraint matrix $A$ and update $Z$ appropriately. The thresholds for the deletion or deflation process to take place can be set as input parameters.

During deflation, we not only fix one variable ($x_{ij}$) to have the value 1, but also fix several other variables to have the value 0. Namely, we fix all variables corresponding to arcs $(i,k)$ for $k \neq j$, $(k,j)$ for $k \neq i$, and $(j,i)$ to have the value 0. Whenever we perform deflation the information about the deflated arcs are stored in order to construct a HC in the original graph once a HC is found in the reduced graph.

After performing deletion or deflation, a reduced vector $\bar{x}$ is obtained, which is  infeasible in the resultant smaller dimension problem. In the stochastic case obtaining a feasible point is trivial since only the variables in the specific rows where fixing has occurred need to adjusted. The simplest way is to multiply the remaining variables in an impacted  rows by $(1/(1-\bar{x}_{ij}))$, where $\bar{x}_{ij}$ is the value of the variable that has been fixed. Note that this increases the remaining variables so will not trigger another deletion in the row. It is possible it triggers a deflation, but this is unlikely. In the doubly stochastic case, recovering feasibility may result in many or all of the variables may be impacted, even for a single variable being deleted. Define
$s := e - A\bar{x}$ to be the error induced by such a process. Note that $s$ is a nonnegative vector
in the case of both deletion and deflation. Then, we find a new $x$ such that $x \in \mathcal{DS}$,
and $|x - \bar{x}| < \varepsilon$, where the size of $\varepsilon$ depends on the deletion or deflation thresholds chosen. The interpretation of $x$ is that it is a point that satisfies the $\mathcal{DS}$ constraints, and is as close as possible to the point we obtained after deleting or deflating.

We determine $x$ by first defining  $v \ge 0$ and $u \ge 0$ so that $x = \bar{x} + v - u$.
Define $x_{min}$ to be a constant which, initially, is equal to the smallest element of $\bar{x}$.
Then, we solve
\begin{eqnarray}&\min\limits_{u,v,\gamma} \rho\gamma + e^T(u + v)& \label{eq-LP_scaling1}\\
&\mbox{s.t.}& \nonumber\\
&Au - Av + \gamma s & = s,\label{eq-LP_scaling2}\\
x_{min}e - \bar{x} \leq & u - v & \leq e - \bar{x},\label{eq-LP_scaling3}\\
0 \leq & u,v,\gamma & \leq 1,\label{eq-LP_scaling4}
\end{eqnarray}

where $\rho$ is chosen large enough that $\gamma$ is reduced to 0 whenever possible. Constraints (\ref{eq-LP_scaling3}) are designed to ensure that $x \in Int(\mathcal{DS})$. However, it may be impossible to satisfy the above constraints for a value of $\gamma = 0$ because some variables may need to be 0 or a value very close to 0. In this case, we reduce $x_{min}$ and solve the LP again, continuing this process until we obtain a solution with $\gamma = 0$.

If $\gamma \ne 0 $ unless we set $x_{ij} = 0$ for some $i$ and $j$, then we delete these variables, as they cannot be nonzero in a Hamiltonian cycle (or in fact any $\mathcal{DS}$ point) containing the currently fixed arcs.

\section*{Rounding}

At each iteration we test if a HC can be obtained by a simple rounding procedure. In the $\mathcal{S}$ case we set the largest element in each row to one, starting with the row with the largest overall element. If the largest element happens to already have a unit element in the column  we set the second largest in that row to one and so on. A similar procedure is used in the $\mathcal{DS}$ case except in this case setting an element to unity induces more elements being set to zero. Moreover, we now fail to satisfy the $\mathcal{DS}$ constraints and so rebalancing is not done.

Obviously, we could use more sophisticated rounding methods, which may allow us to identify a HC earlier. One potential improvement of this method would be to solve a heuristic at the completion of each iteration, using the current point ${\bf x}$, that tries to find a nearby HC. Such a hybrid approach was considered by Eshragh et al. \cite{Ali}, with promising results. This has not been explored since we are interested in testing our algorithm without taking specific advantage of features possibly not common to other problems.

Below we outline the structure of the algorithm, which we term DIPA (Determinant Interior Point Algorithm).

\vspace*{-0.25cm}
{\scriptsize\begin{center}\begin{tabular}{|l|}\hline
\\
\hspace*{0.4cm}{\bf Input}: $\Gamma$. \;\;\;\;\; {\bf Output}: HC found or not found.\\
\\
\hspace*{0.4cm}{\bf begin}\\
\hspace*{1.1cm}${\bf x} \leftarrow $ Find initial interior point\\
\hspace*{1.1cm}${\bf x} \leftarrow $ Find barrier point\\
\hspace*{1.1cm}$\mu \leftarrow \mu_{initial}$\\
\hspace*{1.1cm}{\bf while HC has not been found}\\
\hspace*{1.8cm}{\bf if} {\bf x} is a local min of barrier function\\
\hspace*{2.5cm}Reduce $\mu$\\
\hspace*{2.5cm}{\bf if} $\mu$ too small, {\bf return} no HC found, converged to non-HC local min \;\;\;\;\; \\
\hspace*{1.8cm}{\bf else}\\
\hspace*{2.5cm}{\bf if} reduced Hessian positive definite\\
\hspace*{3.2cm}${\bf d} \leftarrow $ Find descent direction\\
\hspace*{2.5cm}{\bf else}\\
\hspace*{3.2cm}${\bf d} \leftarrow $ Find direction of negative curvature\\
\hspace*{2.5cm}{\bf end}\\
\hspace*{2.5cm}$\alpha \leftarrow $ Choose step length\\
\hspace*{2.5cm}${\bf x} \leftarrow {\bf x} + \alpha d$\\
\hspace*{2.5cm}{\bf begin} deflation/deletion check\\
\hspace*{3.2cm}{\bf if} any variables are above deflation threshold\\
\hspace*{3.9cm}Perform deflation\\
\hspace*{3.2cm}{\bf else if} any variables are below deletion threshold\\
\hspace*{3.9cm}Perform deletion\\
\hspace*{3.2cm}{\bf end}\\
\hspace*{3.2cm}{\bf if} any deletion or deflation occured\\
\hspace*{3.9cm}{\bf if} the underlying graph is no longer connected, {\bf return} no HC found\\
\hspace*{3.9cm}{\bf else}\\
\hspace*{4.6cm}$x \leftarrow $ nearby feasible point\\
\hspace*{4.6cm}{\bf Repeat} deflation/deletion step\\
\hspace*{3.9cm}{\bf end}\\
\hspace*{3.2cm}{\bf end}\\
\hspace*{2.5cm}{\bf end}\\
\hspace*{1.8cm}{\bf end}\\
\hspace*{1.8cm}{\bf if} {\bf x} rounded to HC, {\bf return} HC\\
\hspace*{1.1cm}{\bf end}\\
\hspace*{0.4cm}{\bf end}\\
\\
\hline\end{tabular}\\
\vspace*{0.25cm}{\bf DIPA algorithm.}\end{center}}

\section{Numerical experiments}

In order to investigate the character and to test the performance of the algorithm, we generated a test set of 500 problem instances. Specifically, we randomly generated 50 problems for each of 10, 20, ..., and 100 nodes with node degree between 3 and 6. The computer used for performing all the experiments was a PC with Intel\circledR \;Core$^{\mbox{\scriptsize{TM}}}$ i7-4600U CPU, 2.70 GHz, 16 GB of RAM, and running on the operating system Windows 8.1 Enterprise. DIPA was implemented in MATLAB R2014b, with all LPs solved by IBM ILOG CPLEX Optimization Studio v12.6 via its Concert interface to MATLAB.

\begin{figure}[h]\hspace*{-1.93cm}\includegraphics[scale=0.425]{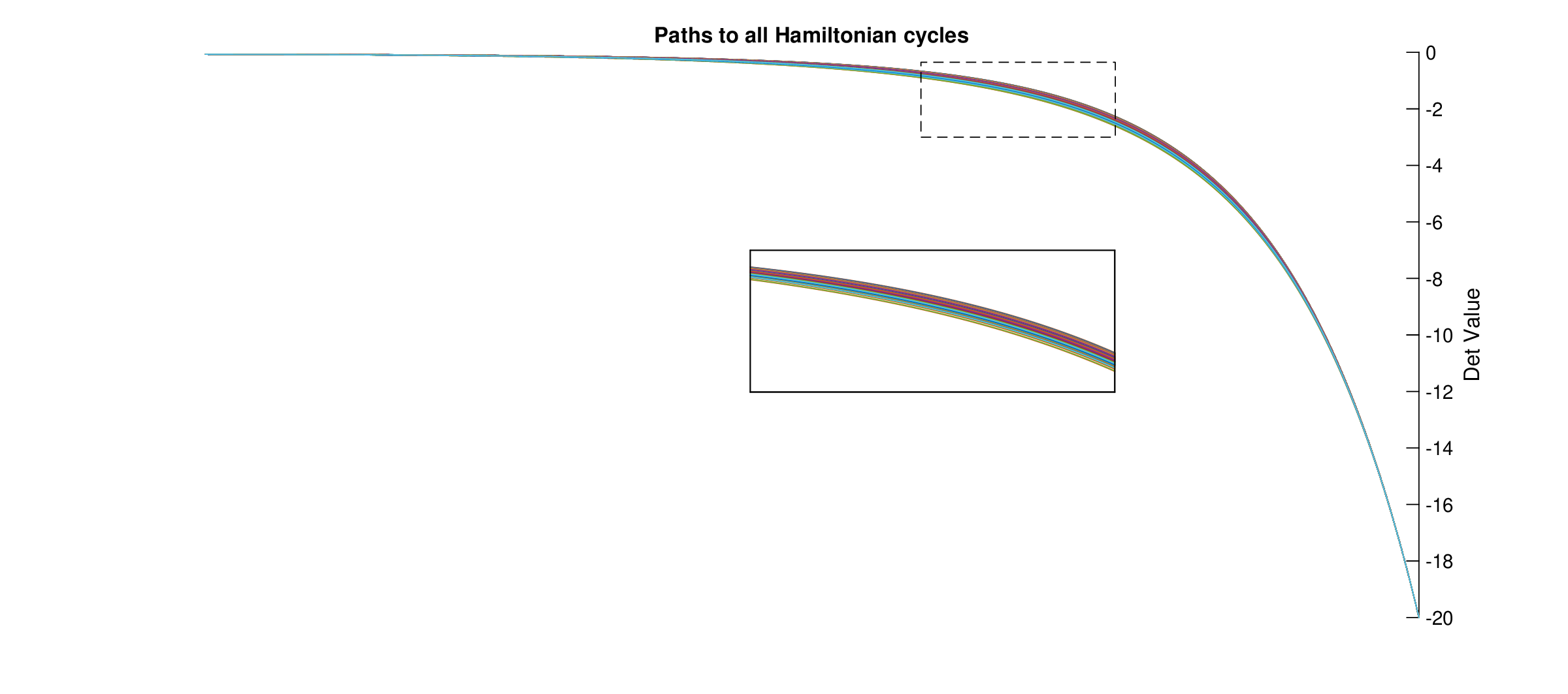}
\caption{Paths showing how the determinant function behaves as it goes, in a straight line, from the initial point to all Hamiltonian cycles of a 20-node graph.}\label{fig-paths}\end{figure}

The choice of initial $\mu$ and the rate of reduction of $\mu$ did not prove to be critical. For successful runs once $\mu$ had been reduced to a sufficient level for the reduced Hessian to be indefinite it almost always remained indefinite until it stopped. Hence, in these cases no further reduction in $\mu$ was required. In Figure \ref{fig-paths} we show how the determinant function behaves as it goes, in a straight line, from the initial point to all of the HCs of a 20 node graph. It can be seen that all curves are remarkably similar, and that each can be reached by going down a direction of negative curvature. Also the degree of curvature increases the closer the point gets to the HC.

We attempted to solve the 500 problems in the test set with the algorithm applied to the stochastic and doubly stochastic form of the problem. In both cases an attempt was made to perform neither a deletion or deflation by setting extreme values for the relevant parameters that invoke these steps within the algorithm.

\begin{table}[h]
\begin{center}
\hspace*{-0.7cm}\begin{tabular}{|c|c|c|c|c|c|c|c|c|c|c|}\hline
\rowcolor[gray]{1}$N$ & {\bf 10} & {\bf 20} & {\bf 30} & {\bf 40} & {\bf 50} & {\bf 60} & {\bf 70} & {\bf 80} & {\bf 90}  & {\bf 100} \\
\rowcolor[gray]{0.9}\hline Solved $\mathcal{S}$ & 47 & 26 & 19 & 13 & 9 & 10 & 7 & 9 & 6 & 4\\
\rowcolor[gray]{1.0}\hline Solved $\mathcal{DS}$ & 50 & 50 & 47 & 48 & 46 & 39 & 45 & 45 & 43 & 41\\
 \hline\end{tabular}
\end{center}
\caption{Numbers of graphs solved when deletions and deflations are suppressed for the stochastic and doubly-stochastic forms of the problem, for graphs of order $N$.}
\label{tab-results}\end{table}

The results given in Table \ref{tab-results} are unambiguous. It is clear that the doubly stochastic form of the problem is far superior. We did do further tests on the stochastic case by varying the adjustable parameters, but the gap in performance was far too large to bridge, and so the remaining tables of results for the DIPA algorithm are for the $\mathcal{DS}$ constrained problem.

A common technique in global optimization when a global minimizer is not found is to try random sets of initial points. Since our algorithm requires a neutral initial point, this is not an option. Instead, we discovered that altering the adjustable parameters and options induced a large variance in which problems were successfully solved.
Since not all problems were solved, we varied some of the options and adjustable parameters. We also introduced two other options. In addition to the use of an LP routine to solve (\ref{eq-LP_scaling1})--(\ref{eq-LP_scaling4}) in order to obtain an interior point after deletion or deflation, we could instead use the CPLEX QP routine to find an interior feasible point of minimum Euclidean length from the current iterate. Since the step from the initial point to the solution of both the QP and LP is very small, one would not expect this to have any impact on efficiency. Indeed, the infeasibility in the linear equalities is very tiny, and this needs to be handled with care since CPLEX can simply treat it as being sufficiently small to be a solution. The other option was to remove one of the variables altogether. Doing so does not alter whether or not a graph has a HC since the reverse cycle exists. One would not expect removing one variable to have a measurable difference in performance. We were interested whether this would impact the set of graphs solved.

We ran the algorithm for several different settings of parameters and options. The algorithm typically ran just as efficiently for any alternative setting of the parameters and options chosen. Graphs that are not successfully solved were attempted again with a different set of parameters. With the option to vary the use of the LP or QP, or different deflation of deletion choices, the attempt to obtain a solution can be made by restarting at the point this option is first used.

The adjustable input parameters and options we considered are:

\begin{itemize}\item Initial $\mu$
\item Amount by which $\mu$ is reduced when a local minimum is reached
\item $\alpha$ - proportion of largest feasible step along the search direction
\item Choose to use a LP or QP to find a feasible point after deletion or deflation
\item Deflation threshold
\item Deletion threshold
\item Whether to use barrier function on $ x \le e$
\item Whether to remove one variable from the problem
\end{itemize}

The choice of $\alpha$ did have an impact on which graphs were solved, but setting $\alpha$ significantly less than one meant that the algorithm usually took longer to converge. In all results reported we set the initial $\mu = 0.01$, the reduction multiplier on $\mu$ to be 0.1, and $\alpha = 0.9$. It appears that deflation is more useful than deletion, so we set the deletion threshold to be very tiny, at 0.00001. In Table \ref{tab-results2} we report the number of graphs from each test set that are solved for four settings. Specifically, we choose between obtaining an interior point via the LP or the QP,  and we set the deflation threshold to be either 0.9 or 0.95. In Table \ref{tab-results3} we also report the number of graphs that are solved after combinations of two of these approaches are employed, and then the final column after all four are employed. In all cases, we include the barrier function on $x \le e$ and remove one variable from each problem.

\begin{table}[h!]
\begin{center}
\hspace*{-0.7cm}\begin{tabular}{|c|c|c|c|c|}\hline
\rowcolor[gray]{1} & Solved & Solved & Solved & Solved\\
\rowcolor[gray]{1}$N$ & (LP, def=0.9) & (LP, def=0.95) & (QP, def=0.9) & (QP, def=0.95)\\
\hline \rowcolor[gray]{0.9}{\bf 10} & 50 & 50 & 50 & 50 \\
\hline \rowcolor[gray]{1.0}{\bf 20} & 49 & 50 & 50 & 49 \\
\hline \rowcolor[gray]{0.9}{\bf 30} & 50 & 49 & 47 & 46 \\
\hline \rowcolor[gray]{1.0}{\bf 40} & 45 & 46 & 42 & 43 \\
\hline \rowcolor[gray]{0.9}{\bf 50} & 42 & 41 & 46 & 45 \\
\hline \rowcolor[gray]{1.0}{\bf 60} & 45 & 39 & 39 & 40 \\
\hline \rowcolor[gray]{0.9}{\bf 70} & 40 & 36 & 40 & 35 \\
\hline \rowcolor[gray]{1.0}{\bf 80} & 40 & 39 & 32 & 36 \\
\hline \rowcolor[gray]{0.9}{\bf 90} & 29 & 32 & 40 & 36 \\
\hline \rowcolor[gray]{1.0}{\bf 100} & 34 & 31 & 28 & 35 \\
\hline\end{tabular}
\end{center}
\caption{Numbers of graphs solved from each set of 50 graphs of order $N$ for different choices of parameters.}
\label{tab-results2}\end{table}

\begin{table}[h!]
\begin{center}
\hspace*{-0.7cm}\begin{tabular}{|c|c|c|c|c|c|}\hline
\rowcolor[gray]{1} & LP, def 0.9 & QP, def=0.9 & LP, def=0.9 & LP, def=0.95 & Combined\\
\rowcolor[gray]{1}$N$ & LP, def=0.95 & QP, def=0.95 & QP, def=0.9 & QP, def=0.95 & (all four)\\
\rowcolor[gray]{0.9}\hline {\bf 10} & 50 & 50 & 50 & 50 & 50 \\
\rowcolor[gray]{1.0}\hline {\bf 20} & 50 & 50 & 50 & 50 & 50 \\
\rowcolor[gray]{0.9}\hline {\bf 30} & 50 & 49 & 50 & 50 & 50 \\
\rowcolor[gray]{1.0}\hline {\bf 40} & 49 & 49 & 49 & 50 & 50 \\
\rowcolor[gray]{0.9}\hline {\bf 50} & 50 & 50 & 50 & 49 & 50 \\
\rowcolor[gray]{1.0}\hline {\bf 60} & 48 & 49 & 50 & 49 & 50 \\
\rowcolor[gray]{0.9}\hline {\bf 70} & 46 & 46 & 46 & 44 & 49 \\
\rowcolor[gray]{1.0}\hline {\bf 80} & 48 & 43 & 47 & 47 & 50 \\
\rowcolor[gray]{0.9}\hline {\bf 90} & 43 & 42 & 46 & 48 & 49 \\
\rowcolor[gray]{1.0}\hline {\bf 100} & 43 & 45 & 44 & 43 & 48 \\
\hline\end{tabular}
\end{center}
\caption{Spread of graphs of order $N$ solved over multiple runs with different choices of parameters.}
\label{tab-results3}\end{table}

In Table \ref{tab-results4} we report similar results for the case when no deflation is performed. Under such circumstances the option of whether to use an LP or QP is irrelevant. Instead, we vary whether or not to have the barrier term on the upper bound ($x \le e$) and whether or not to remove one variable. Table \ref{tab-results5} reports the results of the combinations.

\begin{table}[h!]
\begin{center}
\hspace*{-0.7cm}\begin{tabular}{|c|c|c|c|c|}\hline
\rowcolor[gray]{1} & No Upper Log & No Upper Log & Upper Log & Upper Log\\
\rowcolor[gray]{1}$N$ & No Removed Var & Removed Var & No Removed Var & Removed Var\\
\rowcolor[gray]{0.9}\hline {\bf 10} & 50 & 50 & 50 & 50 \\
\rowcolor[gray]{1.0}\hline {\bf 20} & 50 & 49 & 50 & 50 \\
\rowcolor[gray]{0.9}\hline {\bf 30} & 49 & 49 & 49 & 47 \\
\rowcolor[gray]{1.0}\hline {\bf 40} & 48 & 45 & 48 & 48 \\
\rowcolor[gray]{0.9}\hline {\bf 50} & 49 & 47 & 47 & 46 \\
\rowcolor[gray]{1.0}\hline {\bf 60} & 43 & 45 & 46 & 39 \\
\rowcolor[gray]{0.9}\hline {\bf 70} & 47 & 44 & 44 & 45 \\
\rowcolor[gray]{1.0}\hline {\bf 80} & 45 & 44 & 41 & 45 \\
\rowcolor[gray]{0.9}\hline {\bf 90} & 41 & 43 & 44 & 43 \\
\rowcolor[gray]{1.0}\hline {\bf 100} & 42 & 40 & 42 & 41 \\
\hline\end{tabular}
\end{center}
\caption{Numbers of graphs solved from each set of 50 graphs of order $N$ for different choices of parameters, with deflation and deletion prevented.}
\label{tab-results4}\end{table}

\begin{table}[h!]
\begin{center}
\hspace*{-0.7cm}\begin{tabular}{|c|c|c|c|c|c|}\hline
\rowcolor[gray]{1} & Upper Log  & Removed Var & No Upper Log & No Removed Var & Combined\\
\rowcolor[gray]{1}$N$ & Combo & Combo & Combo & Combo & (all four)\\
\rowcolor[gray]{0.9}\hline {\bf 10} & 50 & 50 & 50 & 50 & 50 \\
\rowcolor[gray]{1.0}\hline {\bf 20} & 50 & 50 & 50 & 50 & 50 \\
\rowcolor[gray]{0.9}\hline {\bf 30} & 50 & 50 & 50 & 49 & 50 \\
\rowcolor[gray]{1.0}\hline {\bf 40} & 50 & 49 & 49 & 48 & 50 \\
\rowcolor[gray]{0.9}\hline {\bf 50} & 50 & 49 & 50 & 49 & 50 \\
\rowcolor[gray]{1.0}\hline {\bf 60} & 49 & 48 & 50 & 48 & 50 \\
\rowcolor[gray]{0.9}\hline {\bf 70} & 49 & 50 & 50 & 48 & 50 \\
\rowcolor[gray]{1.0}\hline {\bf 80} & 48 & 50 & 49 & 47 & 50 \\
\rowcolor[gray]{0.9}\hline {\bf 90} & 49 & 50 & 50 & 50 & 50 \\
\rowcolor[gray]{1.0}\hline {\bf 100} & 47 & 47 & 48 & 46 & 50 \\
\hline\end{tabular}
\end{center}
\caption{Spread of graphs of order $N$ solved over multiple runs with different choices of parameters, with deflation and deletion prevented.}
\label{tab-results5}\end{table}

\begin{figure}[h]\hspace*{-1.93cm}\includegraphics[scale=0.425]{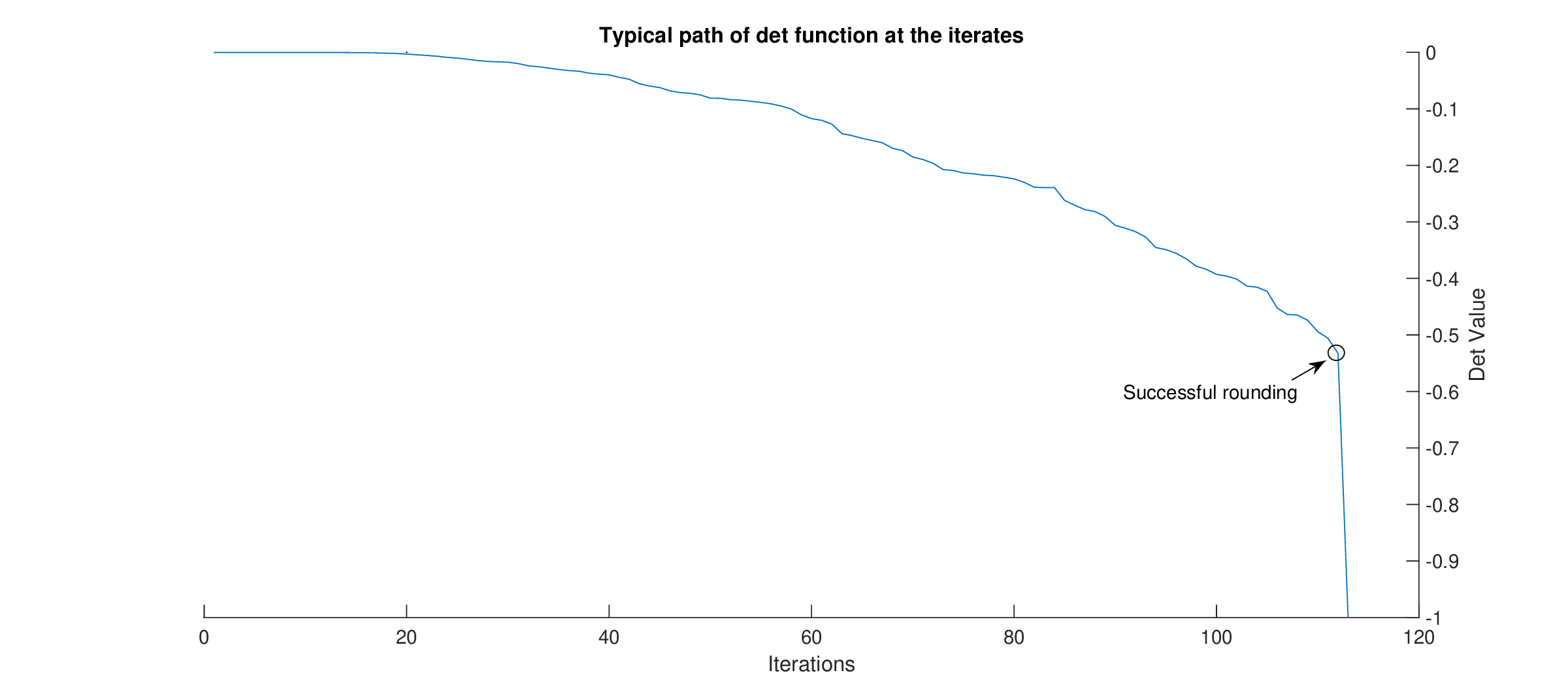}
\caption{A typical path of the determinant function at the iterates of the algorithm. After 112 iterations the rounding procedure finds a HC.}\label{fig-descent}\end{figure}

Not using deflation typically increased the number of iterations about 50\% over using deflation. In Figure \ref{fig-descent} we show how our algorithm typically converges when a HC is found. Since we are converging to a vertex the nature of converges differs considerably from a typical minimization algorithm where the reduction made in the objective slows as the solution is approached.

\section*{Results for the Solver BONMIN}
To assess the difficulty of the problems, we also tried to solve them using BONMIN \cite{bonmin1,bonmin2}, which is designed both to solve nonlinear problems with continuous variables and problems with discrete variables.  Specifically, we used the implementation of BONMIN included in the OPTI Toolbox v2.16 \cite{opti}. It is one of the few codes in the public domain. Note that BONMIN is only designed to find a local minimizer. Since BONMIN does not require a neutral point, we used four initial points to improve the chance of finding a global minimizer.  For each of the 500 graphs, we first generated a feasible point, and then obtained four feasible points by computing four random feasible directions and taking a partial step to the boundary. BONMIN has many options and we used the default options with the exception in a few cases where the time limit to solve the problem was increased from 1,000 secs to 10,000 secs. Whenever BONMIN terminated, we applied the same rounding algorithm described above to check if a Hamiltonian cycle had been found.


In Table \ref{tab-results6} we report the number of problems for which a global minimum was successfully found, for each of the four starting points. In Table \ref{tab-results7}, we report the number of graphs successfully solved for four different combinations of runs, and then the combined set of graphs solved over all runs.

\begin{table}[h!]
\begin{center}
\hspace*{-0.7cm}\begin{tabular}{|c|c|c|c|c|}\hline
\rowcolor[gray]{1} & Starting & Starting & Starting & Starting \\
\rowcolor[gray]{1}$N$ & Point 1 & Point 2 & Point 3 & Point 4\\
\rowcolor[gray]{0.9}\hline {\bf 10} & 39 & 31 & 36 & 37\\
\rowcolor[gray]{1.0}\hline {\bf 20} & 28 & 30 & 22 & 29\\
\rowcolor[gray]{0.9}\hline {\bf 30} & 24 & 30 & 23 & 26\\
\rowcolor[gray]{1.0}\hline {\bf 40} & 28 & 13 & 25 & 25\\
\rowcolor[gray]{0.9}\hline {\bf 50} & 16 & 18 & 26 & 23\\
\rowcolor[gray]{1.0}\hline {\bf 60} & 23 & 21 & 20 & 21\\
\rowcolor[gray]{0.9}\hline {\bf 70} & 21 & 23 & 14 & 16\\
\rowcolor[gray]{1.0}\hline {\bf 80} & 15 & 10 & 15 & 17\\
\rowcolor[gray]{0.9}\hline {\bf 90} & 6 & 6 & 3 & 6\\
\rowcolor[gray]{1.0}\hline {\bf 100} & 0 & 0 & 0 & 0\\
\hline\end{tabular}
\end{center}
\caption{Spread of graphs of order $N$ solved by BONMIN, with four different starting points.}
\label{tab-results6}\end{table}

\begin{table}[h!]
\begin{center}
\hspace*{-0.7cm}\begin{tabular}{|c|c|c|c|c|c|}\hline
\rowcolor[gray]{1} & First and  & First and & Second and & Third and & Combined\\
\rowcolor[gray]{1}$N$ & Second Pts & Third Pts & Fourth Pts & Fourth Pts & (all four)\\
\rowcolor[gray]{0.9}\hline {\bf 10} & 46 & 45 & 42 & 44 & 47 \\
\rowcolor[gray]{1.0}\hline {\bf 20} & 41 & 33 & 41 & 34 & 47 \\
\rowcolor[gray]{0.9}\hline {\bf 30} & 36 & 37 & 40 & 38 & 48 \\
\rowcolor[gray]{1.0}\hline {\bf 40} & 35 & 40 & 32 & 38 & 47 \\
\rowcolor[gray]{0.9}\hline {\bf 50} & 27 & 34 & 32 & 36 & 43 \\
\rowcolor[gray]{1.0}\hline {\bf 60} & 35 & 32 & 34 & 31 & 46 \\
\rowcolor[gray]{0.9}\hline {\bf 70} & 35 & 29 & 31 & 25 & 42 \\
\rowcolor[gray]{1.0}\hline {\bf 80} & 21 & 26 & 21 & 27 & 36 \\
\rowcolor[gray]{0.9}\hline {\bf 90} & 11 & 9 & 11 & 7 & 16 \\
\rowcolor[gray]{1.0}\hline {\bf 100} & 0 & 0 & 0 & 0 & 0 \\
\hline\end{tabular}
\end{center}
\caption{Spread of graphs of order $N$ solved over multiple runs with BONMIN.}
\label{tab-results7}\end{table}

BONMIN also permits the user to specify that a variable should be binary in the solution. Since BONMIN first solves the relaxed problem, clearly the binary choice will not be worse. Of course, specifying the binary option adds to the time to find a solution, and in our experience this additional time was many times more than the solve time for the relaxed problem. It is possible that the better option is not to specify the binary requirement, but to solve for more initial points instead. In Tables \ref{tab-results8} and \ref{tab-results9} we display the equivalent results of the previous two tables once binary solutions are requested. It may be observed from Tables \ref{tab-results7} and \ref{tab-results8} that the combined result of four solves of the relaxed problem is often better than a single solve with the binary requirement.

\begin{table}[h!]
\begin{center}
\hspace*{-0.7cm}\begin{tabular}{|c|c|c|c|c|}\hline
\rowcolor[gray]{1} & Starting & Starting & Starting & Starting \\ \rowcolor[gray]{1}$N$ & Point 1 & Point 2 & Point 3 & Point 4\\ \rowcolor[gray]{0.9}\hline {\bf 10} & 49 & 47 & 46 & 47\\ \rowcolor[gray]{1.0}\hline {\bf 20} & 44 & 41 & 44 & 45\\ \rowcolor[gray]{0.9}\hline {\bf 30} & 41 & 45 & 39 & 46\\ \rowcolor[gray]{1.0}\hline {\bf 40} & 44 & 36 & 39 & 42\\ \rowcolor[gray]{0.9}\hline {\bf 50} & 35 & 32 & 38 & 35\\ \rowcolor[gray]{1.0}\hline {\bf 60} & 31 & 31 & 37 & 35\\ \rowcolor[gray]{0.9}\hline {\bf 70} & 29 & 33 & 19 & 29\\ \rowcolor[gray]{1.0}\hline {\bf 80} & 27 & 26 & 30 & 32\\ \rowcolor[gray]{0.9}\hline {\bf 90} & 18 & 19 & 20 & 17\\ \rowcolor[gray]{1.0}\hline {\bf 100} & 8 & 9 & 8 & 5\\ \hline\end{tabular} \end{center} \caption{Spread of graphs of order $N$ solved by BONMIN with binary variables, with four different starting points.} \label{tab-results8}\end{table}

\begin{table}[h!]
\begin{center}
\hspace*{-0.7cm}\begin{tabular}{|c|c|c|c|c|c|}\hline
\rowcolor[gray]{1} & First and  & First and & Second and & Third and & Combined\\ \rowcolor[gray]{1}$N$ & Second Pts & Third Pts & Fourth Pts & Fourth Pts & (all four)\\ \rowcolor[gray]{0.9}\hline {\bf 10} & 50 & 50 & 50 & 50 & 50 \\ \rowcolor[gray]{1.0}\hline {\bf 20} & 49 & 48 & 49 & 49 & 50 \\ \rowcolor[gray]{0.9}\hline {\bf 30} & 49 & 50 & 50 & 50 & 50 \\ \rowcolor[gray]{1.0}\hline {\bf 40} & 49 & 48 & 46 & 48 & 50 \\ \rowcolor[gray]{0.9}\hline {\bf 50} & 42 & 44 & 45 & 47 & 50 \\ \rowcolor[gray]{1.0}\hline {\bf 60} & 43 & 44 & 46 & 42 & 49 \\ \rowcolor[gray]{0.9}\hline {\bf 70} & 46 & 38 & 37 & 42 & 48 \\ \rowcolor[gray]{1.0}\hline {\bf 80} & 40 & 41 & 42 & 38 & 48 \\ \rowcolor[gray]{0.9}\hline {\bf 90} & 28 & 34 & 31 & 26 & 40 \\ \rowcolor[gray]{1.0}\hline {\bf 100} & 14 & 14 & 14 & 13 & 19 \\ \hline\end{tabular} \end{center} \caption{Spread of graphs of order $N$ solved over multiple runs with BONMIN with binary variables.} \label{tab-results9}\end{table}

Finally, in Table \ref{tab-results10} we display the success rate of finding global minimizers when only $\mathcal{S}$ constraints are used. In this experiment, we also requested binary solutions. It is clear from Table \ref{tab-results10} that further experimentation with BONMIN on the stochastic form of the problem was unnecessary.

\begin{table}[h!]
\begin{center}
\hspace*{-0.7cm}\begin{tabular}{|c|c|c|c|c|c|}\hline
\rowcolor[gray]{1} & Starting & Starting & Starting & Starting & Combined\\
\rowcolor[gray]{1}$N$ & Point 1 & Point 2 & Point 3 & Point 4 & Solves\\
\rowcolor[gray]{0.9}\hline {\bf 10} & 15 & 25 & 20 & 22 & 45\\
\rowcolor[gray]{1.0}\hline {\bf 20} & 0 & 1 & 1 & 0 & 2\\
\rowcolor[gray]{0.9}\hline {\bf 30} & 0 & 0 & 0 & 0 & 0\\
\rowcolor[gray]{1.0}\hline {\bf 40} & 0 & 0 & 0 & 0 & 0\\
\rowcolor[gray]{0.9}\hline {\bf 50} & 0 & 0 & 0 & 0 & 0\\
\hline\end{tabular}
\end{center}
\caption{Spread of graphs of order $N$ solved by BONMIN with binary variables and $\mathcal{S}$ constraints, with four different starting points, and the combined graphs solved over the four separate runs.}
\label{tab-results10}\end{table}

\section{Concluding remarks}

The performance of BONMIN on this class shows that finding a global minimizer is a challenge, even when binary variables are specified. Also, the results for  the problem with $\mathcal{DS}$ constraints are far better than the problem with $\mathcal{S}$ constraints. This supports the view that different forms of a problem with identical complexity properties can have quite differing performance. It was not the intent of this paper to show or suggest that our algorithm is competitive with alternative algorithms to find a Hamiltonian cycle. However, it is quite distinct from other methods and there is much that can be done to improve its performance. More constraints can be added. For example, $P^T P = I$, and we know that twin variables, say $x_i$ and $x_j$ must satisfy $1 - x_i - x_j \ge 0$ and $x_i x_j = 0$. The product from the latter constraints can be added directly to the objective. However, we plan to try and find a form of the problem that eliminates the occurrence of reverse cycles and so eliminates twin variables from the problem. The new variables would be the elements of $Q$, where $Q = 0.5(P + P^T)$. Knowing $Q$ it is trivial to find the elements of $P$. The current problem has a dense Hessian matrix. Although we have shown how all the elements can be computed efficiently, it still leaves a dense matrix, which has computational implications when computing the search direction for large problems. We are investigating transformations that should lead to the Hessian being sparse. Also, if the conjugate gradient algorithm is used to compute the search direction and direction of negative curvature, it may be possible to compute $H v$ efficiently even when $H$ is dense.

Although we have addressed the Hamiltonian cycle problem, an equally important interest is developing algorithms to determine global minimizers for other relaxed discrete problems. Many of the issues that arise in such problems are identical to those arising in the HC problem. For example, lots of global minimizers and hence lots of stationary points that have reduced Hessians that are almost positive definite (one negative eigenvalue). Moreover, symmetry is also present. Problems such as the frequency assignment problem have an equally good solution simply from any permutation of a known solution. Also, for constrained problems of this type, solutions are typically at a highly degenerate vertex. We are encouraged by the success of the algorithm we have developed, which has demonstrated the ability to find global minimizers of highly nonlinear and  nonconvex problems with several hundred binary variables. As already noted a common technique used in algorithms to solve global optimization problems is to use multiple starting points. The approach we advocate requires a neutral starting point, which is usually unique. We have demonstrated that an equally good alternative (and possibly better) is to vary some of the parameters and options that algorithms to solve such problems typically have.  We have shown that very small changes both to the strategy and flexible parameters leads to distinct solution, enabling us to reduce significantly the number of problems on which we fail to find a global minimizer. Moreover, these variations do not lead to less efficient methods. A possible explanation is that the algorithm generates iterates that lie on the edges of two domains that contain different global minimizers.

We think the form of the Hamiltonian cycle problem we present makes an excellent test problem both for algorithm designed to solve discrete problems and algorithms to find global minimizers. The problems are simple to generate and the degree of difficult is a function of the number of arcs. It is easy to test whether the best solution has been found without compromising how the test problem is generated. Another point of interest is whether transforming a discrete problem into a continuous problem whose extreme discrete points are minimizers (implying rounding is always an improvement) is always worthwhile. What seems more fruitful is adding more constraints.

\section*{Data Availability}

All data generated or analysed during this study are included in this published article and its supplementary information files.

\end{document}